\documentclass[12pt,twoside,makeidx]{amsart}
\usepackage{latexsym}
\usepackage{amsmath}
\usepackage{epsfig}
\usepackage{epic}
\usepackage{amssymb}
\usepackage{enumerate, afterpage, float, tikz, url, mparhack}

\usepackage{hyperref}

\newtheorem{theorem}{Theorem}[section]
\newtheorem{lemma}[theorem]{Lemma}
\newtheorem{corollary}[theorem]{Corollary}
\newtheorem{proposition}[theorem]{Proposition}

\newtheorem{sublemma}{}[theorem]

\theoremstyle{definition}

\theoremstyle{remark}

\numberwithin{equation}{section}

\newcommand{\ba}{\backslash}

\newcommand{\co}{{\rm co}}
\newcommand{\si}{{\rm si}}
\newcommand{\cl}{{\rm cl}}

\newcommand{\subproof}{\begin{proof}[Subproof]}

\newcommand{\del}{\backslash}

\newcommand{\bF}{\mathbb F}

\newcommand{\cT}{\mathcal{T}}

\begin{document}

\sloppy

\title{Certifying Non-representability of Matroids Over
Prime Fields}

\author{Jim Geelen}
\address{Department of Combinatorics and Optimization\\
University of Waterloo\\
Waterloo, Ontario, Canada}
\email{jfgeelen@math.uwaterloo.ca}
\thanks{The first author was supported by a grant from the
National Sciences and Engineering Research Council of Canada.}


\author{Geoff Whittle}
\address{School of Mathematics, Statistics and Operations Research\\
Victoria University of Wellington,  New Zealand}
\email{geoff.whittle@vuw.ac.nz}
\thanks{The second author was supported by a grant from the Marsden
Fund of New Zealand.}

\date{}

\begin{abstract}
It is proved that, for a prime number $p$, showing that an
$n$-element matroid is not representable over $GF(p)$
requires only $O(n^2)$ rank evaluations.
\end{abstract}

\maketitle



\section{Introduction}

We prove the following theorem.

\begin{theorem}
\label{certify}
Let $p$ be a prime number. Then proving that
an $n$-element matroid is not 
$GF(p)$-representable requires at most $O(n^2)$ 
rank evaluations.
\end{theorem}

Theorem~\ref{certify} is a corollary of a more technical 
result.
This technical result---stated as 
Theorem~\ref{sufficient}---gives sufficient
conditions for a matroid to have a certificate of
non-representability over a finite field using at most $O(n^2)$
rank evaluations. 
That these sufficient conditions are satisfied
for prime fields follows from results in two
other papers. 
The first is a theorem of Geelen and Whittle \cite{gewh1}
that shows that matroids satisfying a certain strengthening
of $3$-connectivity have a bounded number of inequivalent
$GF(p)$-representations. The second is a theorem of
Ben David and Geelen  \cite{dage} that shows
that excluded minors for $GF(p)$-representability
cannot have arbitrarily long nested
sequences of separations of bounded
order, in particular they cannot have long 
nested sequences of 3-separations. 
Theorem~\ref{sufficient} is stated at a higher level
of generality than we need
for this paper. It is hoped that by doing this it may 
prove to be a useful tool
in a future extension of Theorem~\ref{certify} to the non-prime
case.

Before proceeding to the technicalities we give some background to
the problems addressed by this paper. In general it is not
possible to describe matroids with an input that is polynomial
in the size of the ground set---there are far too many matroids.
Because of this, it is common to consider oracles, typically rank
oracles. The oracle has a matroid $M$ in mind and we can
ask the oracle for the rank of any subset that we care about.
A call to the oracle counts as a single step in an algorithm.
Seymour \cite{seymour1} showed that, in general, it requires exponentially
many calls to a rank oracle to decide if a matroid is binary.
Seymour's techniques can be generalised to obtain
an analogous negative result for any other field. While this
extension is widely known in the matroid community, as far as we know
it does not appear in print and we give the details in 
Section~\ref{cert-rep}. 

The situation for certifying non-representability is much more
hopeful. Consider binary matroids. We know \cite{tutte1}
that $U_{2,4}$ is the unique excluded minor for binary matroids.
It is also the case that a matroid can be proved to have a 
$U_{2,4}$-minor using eight rank evaluations. Therefore
a matroid can be proved to be non-binary using 
only $O(1)$ rank evaluations. Rota \cite{rota1} conjectured that
for each finite field ${\mathbb F}$, there are, up to isomorphism, only 
finitely many excluded minors for the class of 
${\mathbb F}$-representable matroids.
If Rota's Conjecture were true, then it would require only
$O(1)$ rank evaluations to certify non-representability for 
matroids over any
finite field. While evidence is accumulating for the truth
of Rota's Conjecture, to date it is still only known to hold
for $GF(2)$, $GF(3)$ \cite{bixby1, kase1, seymour2} 
and $GF(4)$ \cite{ggk}.

An alternative approach to finding a short 
proof that a matroid is not binary matroid is to 
attempt to build a binary matrix $A$ with the property
that $M$ is binary if and only if $M=M[A]$.
Such a matrix can easily be constructed by choosing a
basis $B=(b_1,b_2,\ldots,b_r)$ and considering 
basic circuits relative to $B$:
the element $a_{ij}$ of $A$ is nonzero if and only 
if $b_i$ is in the basic circuit of the 
element labelling the $i$th column of $A$. 
Moreover, it can be shown that
$M\neq M[A]$ with a single rank evaluation. This is an
approach that has the potential to be extended to other fields,
and does not rely on the truth of Rota's Conjecture. 

The problem with extending this technique to other fields
is that, in general, a matroid may have many inequivalent 
representations over a field. Recall that two 
$\mathbb F$-representations of a matroid
$M$ are equivalent if one can be obtained from the other 
by elementary row operations (adding one row to the other,
adjoining or deleting a row of zeroes, and scaling a row) 
and column-scaling. For very small fields one can control
the number of inequivalent representations of a matroid.
Indeed matroids are uniquely representable over
$GF(2)$ and $GF(3)$.
Moreover Kahn \cite{kahn1} showed that 3-connected matroids
have at most two inequivalent representations over $GF(4)$ (or
are uniquely representable if field automorphisms are allowed).
In \cite{oxvewh1} it is shown that 3-connected matroids
have at most six inequivalent representations over $GF(5)$.
This is a key ingredient in the proof given in \cite{govw}
that  $O(n^2)$ rank evaluations suffice to prove that a matroid
is not $GF(5)$-representable.

Unfortunately for fields with more than five elements,
no bound can be placed on the number of inequivalent
representations of 3-connected matroids. But, as noted earlier,
it follows from results in \cite{gewh1} that the number of inequivalent
representations of a matroid over a fixed prime field $GF(p)$
becomes
bounded if one raises the connectivity somewhat. We also need
to ensure that whatever connectivity we are dealing with
is possessed by  excluded minors for $GF(p)$. This assurance
is a consequence of the previously mentioned result \cite{dage}
where it is shown that excluded minors over $GF(p)$
do not have arbitrarily long nested sequences of $3$-separations.

\section{Certifying GF$(q)$-representability}
\label{cert-rep}

In this section we give 
examples showing that there is no succinct certificate for
representability over any finite field.
Specifically, we show that, for some $\alpha>1$, proving 
$GF(q)$-representability
for $n$-element matroids requires at least $\alpha^n$ rank-values
in the worst case.  These examples are related to results of 
Seymour~\cite{seymour1} who showed that there is no efficient
algorithm to determine whether or not a matroid is binary.

Let $n\geq 3$ be an integer and
let $N$ be a rank-$n$ matroid with ground set
$\{t,a_1,b_1,a_2,b_2,\ldots,a_n,b_n\}$ such that
\begin{itemize}
\item[(i)] $\{t,a_i,b_i\}$ is a triangle for all 
$i\in\{1,2,\ldots,n\}$, and
\item[(ii)] $r(\cup_{j\in J}\{a_j,b_j\})=|J|+1$ for 
every proper subset $J$ of $\{1,2,\ldots,n\}$.
\end{itemize}
Then the matroid $M=N\ba t$ is a rank-$n$ {\em spike}.  
Each pair $\{a_i,b_i\}$ is a {\em leg}
of the spike. 
For distinct $i,j\in\{1,2,\ldots,n\}$, the
set $\{a_i,b_i,a_j,b_j\}$ is a circuit of $M$.
Any other  non-spanning circuit of a spike,
is a transversal of the legs.

Representable spike can be obtained as follows.
Let $\bF$ be a finite field of order $q\ge 3$ and
let $J=\{j_0,j_1,\ldots ,j_n\}$ be a spanning
circuit of PG$(n-1,q)$.
Now, for each $i\in \{1,2,\ldots,n\}$,
let $a_i$ and $b_i$ be distinct points of $PG(n-1,q)\del\{j_{0}\}$
that are on the line spanned by $j_{0}$ and $j_i$.
Now let $M$ be the restriction of PG$(n-1,q)$ to
$\{a_1,b_1,a_2,b_2,\ldots,a_n,b_n\}$.
The next theorem, which is stronger than we need, shows that it really
does require a lot of work to pin down a $GF(q)$-representable
spike.

\begin{theorem}\label{cert1}
Let $q$ be a prime power and let $M$ be a spike with rank $n> q^3$.
To prove that $M$ is $GF(q)$-representable
requires at least $2^{\frac{n}{2}}$ rank values.
\end{theorem}

Let $M$ be a spike with legs $l_1=\{a_1,b_1\},\ldots,l_n=\{a_n,b_n\}$.
Let $\cT(M)$ denote the set of all dependent
transversals of $(l_1,l_2,\ldots,l_n)$. The following claim is straightforward.

\begin{lemma}
\label{triv1}
If $T_1,T_2\in \cT(M)$, then $|T_1-T_2|\neq 1$.
\end{lemma}

\begin{proof}
If the result fails, then, by symmetry, we may assume that
$T_1=\{a_1,\ldots,a_n\}$ and that $T_2=\{a_1,\ldots,a_{n-1},b_n\}$.
Thus, $r_M(T_1\cup T_2)=n-1$. Then, considering the 
$4$-element circuits $l_i\cup l_n$, we see that
$r(M)=n-1$ contradicting the fact that $r(M)=n$.
\end{proof}

The following result is also straightforward; we leave the proof
to the reader.
\begin{lemma}
\label{triv2}
Let $\cT$ be a set of transversals of $(l_1,l_2,\ldots, l_n)$
such that $|T_1-T_2|\neq 1$ for any $T_1,T_2\in \cT$.
Then there exists a unique spike $M$ on legs
$l_1,l_2,\ldots,l_n$ such that $\cT(M) = \cT$.
\end{lemma}

The complexity results in this section rely on the following lemma.
\begin{lemma}\label{diff}
Let $M_1$ and $M_2$ be GF$(q)$-representable spikes with
legs $l_1,\ldots,l_n$ such that
$\cT(M_1) = \cT(M_2)\cup \{T\}$ for
some transversal $T\not\in \cT(M_1)$ of $(l_1,\ldots,l_n)$.
Then $n\le (q-1)^3+1$.
\end{lemma}

\begin{proof}
We may assume that $T=\{a_1,a_2\ldots,a_{n-1},b_n\}$.  
By Lemma~\ref{triv1},
the transversal $\{a_1,a_2,\ldots,a_n\}$ is independent in
$M_2$ and, hence, also in $M_1$.
Consider representations of $M_1$ and $M_2$; we have:
$$ A_1 = 
\bordermatrix{
 & a_1 &  a_2 &  \cdots & a_n & b_1 & b_2 & \cdots & b_n \cr
 & 1 & 0 & \cdots & 0 & 1+\alpha_1^{-1} & 1 & \cdots & 1 \cr
 & 0 & 1 &  & 0 & 1&1+\alpha_2^{-1} &   & 1 \cr
 & \vdots & &\ddots& & \vdots &&\ddots &\cr
 & 0 & 0 &  & 1 & 1&1 &  & 1+\alpha_n^{-1} 
} \mbox{ and}
$$
$$ A_2 = 
\bordermatrix{
 & a_1 &  a_2 &  \cdots & a_n & b_1 & b_2 & \cdots & b_n \cr
 & 1 & 0 & \cdots & 0 & 1+\beta_1^{-1} & 1 & \cdots & 1 \cr
 & 0 & 1 &  & 0 & 1&1+\beta_2^{-1} &   & 1 \cr
 & \vdots & &\ddots& & \vdots &&\ddots &\cr
 & 0 & 0 &  & 1 & 1&1 & & 1+\beta_n^{-1}, 
}
$$
where $\alpha_i\neq 0$ and $\beta_i\neq 0$ 
for each $i\in\{1,2,\ldots,n\}$.
For $S\subseteq \{1,2,\ldots,n\}$ we let $T_S$ denote the transversal 
$\{a_i\, : \, i\not\in S\}\cup \{b_i\, : \, i\in S\}$.
It is easily verified that $T_S$ is dependent in $M_1$ if and only if
$\sum (\alpha_i\, : \, i\in S) = -1$ and that
$T_S$ is dependent in $M_2$ if and only if
$\sum (\beta_i\, : \, i\in S) = -1$.
Suppose, by way of contradiction, that $n>(q-1)^3+1$.  
Note that there
are at most $(q-1)^2$ distinct pairs $(\alpha_i,\beta_i)$. Therefore,
there exist $S\subseteq \{1,2,\ldots,n-1\}$, 
and $\alpha,\beta\in GF(q)-\{0\}$
such that $|S|=q$ and $(\alpha_i,\beta_i)=(\alpha,\beta)$ 
for each $i\in S$.
Now, $\sum( \alpha_i\, :\, i\in S\cup \{n\}) = q\alpha + \alpha_n =
\alpha_n\neq -1$ and $\sum( \beta_i\, :\, i\in S\cup n) 
= q\beta + \beta_n =\beta_n= -1$. 
So $T_{S\cup\{n\}}$ is dependent in $M_2$ but not in $M_1$,
which is a contradiction.
\end{proof}

\begin{proof}[Proof of Theorem~\ref{cert1}.]
Suppose that $M$ has legs $(l_i\, : \, i=1,\ldots,n)$.
Let $T$ be a transversal of $(l_1,l_2,\ldots,l_n)$.
If $T\in\cT(M)$, then, by Lemma~\ref{triv1} and Lemma~\ref{diff},
relaxing $T$ results in a non-GF$(q)$-representable spike.
Similarly, if $|T-T_1|>1$ for all $T_1\in\cT(M)$, then
restricting $T$ to a circuit results in a non-GF$(q)$-representable spike.
Let $\cT'$ denote the set of all transversals $T$ 
of $(l_1,l_2,\ldots,l_n)$
such that $|T-T_1|>1$ for all $T_1\in\cT(M)$.
Thus, there are $|\cT(M)| + |\cT'|$ non-GF$(q)$-representable spikes
that differ in rank from $M$ only on one set.  Moreover, we have
$ (n+1) (|\cT(M)| + |\cT'|) \ge (n+1)|\cT(M)| + |\cT'| \ge 2^n$.
Thus $|\cT(M)| + |\cT'| \ge 2^n/(n+1) \ge 2^{\frac{n}{2}}$.
Therefore, to distinguish $M$ from each of these 
non-GF$(q)$-representable matroids, we need 
at least $2^{\frac{n}{2}}$ rank values.
\end{proof}

\section{Freedom in matroids}

In this section we review basic material on freedom in matroids
and prove some lemmas that we will need for Theorem~\ref{sufficient}.
The treatment of freedom given here follows \cite{govw}.

Let $M$ be a matroid. Elements $e$ and $f$ of $M$ are 
{\em clones} if swapping the labels of $e$ and $f$ is an
automorphism of $M$. A {\em clonal class} of 
$M$ is a maximal set of elements of $M$ every pair
of which are clones. An element $z$ of $M$ is {\em fixed} in 
$M$ if there is no extension of $M$ by an element $z'$ in which
$z$ and $z'$ are independent clones. Dually, an element
$z'$ of $M$ is {\em cofixed} in $M$ if it is fixed in $M^*$.
Note that if $z$ already has a clone, say $x$, and $\{x,z\}$
is independent, then $z$ is not fixed as we can add a new
element freely on the line through $x$ and $z$.

A flat $F$ of $M$ is {\em cyclic} if it is a union of circuits
of $M$. The next result is straightforward.

\begin{proposition}
\label{free-for-all}
Elements $e$ and $f$ of a matroid $M$ are clones if and only
if they are contained in the same set of cyclic flats.
\end{proposition}

Let $e$ and $f$ be elements of $M$. Then $e$ is {\em freer}
than $f$ if every cyclic flat containing $e$ also contains $f$.
Thus $e$ and $f$ are clones if and only if $e$ is freer than
$f$ and $f$ is freer than $e$. The {\em freedom} of an
element $e$ of $M$ is the maximum size of an independent 
clonal class containing $e$ among all extensions of $M$.
This maximum does not exist if and only if $e$ is a coloop of
$M$; in this case the freedom of $e$ is infinity. Observe that
an element if fixed in $M$ if and only if its
freedom is 0 or 1.

The notion of freedom of an element in a matroid was introduced
by Duke \cite{duke}. While his definition is different from that
given here, it is, in fact, equivalent; see \cite[Lemma~2.8]{govw}.
The next lemma is \cite[Theorem~6.2]{duke}. A proof is also
given in \cite[Lemma~2.9]{govw}.

\begin{lemma}
\label{freer-freedom}
Let $a$ and $b$ be elements of the matroids $M$ such that
$a$ is freer than $b$. Then the freedom of $a$ is at least
the freedom of $b$. Moreover, either $a$ and $b$ are clones
or the freedom of $a$ is greater than the freedom of $b$.
\end{lemma}

For elements $e$ and $f$ of a matroid $M$, it is easily
seen that the freedom of $f$ does not decrease when we delete
$e$. Contraction has a slightly more complicated effect
on freedom. The proof of part (ii) of the next lemma is given 
in \cite[Lemma~2.10]{govw}.

\begin{lemma}
\label{minor-freedom}
Let $e$ and $f$ be elements of the matroid $M$ where 
$f$ has freedom $k$ in $M$.
\begin{itemize}
\item[(i)] The freedom of $f$ in $M\ba e$ is at least $k$.
\item[(ii)] The freedom of $f$ in $M/e$ is at least $k-1$.
Moreover, if $f$ has freedom $k-1$ in $M/e$, then 
$f$ is freer than $e$ in $M$.
\end{itemize}
\end{lemma}

The {\em cofreedom} of an element $e$ of $M$ is the freedom of
$e$ in $M^*$. It is easily seen that the cyclic flats
of $M^*$ are the complements of the cyclic flats of $M$. It
follows from this fact that $e$ is freer than $f$ in $M^*$
if and only if $f$ is freer than $e$ in $M$. Note that an
element is cofixed in $M$ if and only if its cofreedom is either
0 or 1. The next
lemma is the dual of Lemma~\ref{minor-freedom}.

\begin{lemma}
\label{cominor-freedom}
Let $e$ and $f$ be elements of the matroid $M$ where 
$f$ has cofreedom $k$ in $M$.
\begin{itemize}
\item[(i)] The cofreedom of $f$ in $M/ e$ is at least $k$.
\item[(ii)] The cofreedom of $f$ in $M\ba e$ is at least $k-1$.
Moreover, if $f$ has cofreedom $k-1$ in $M\ba e$, then 
$e$ is freer than $f$ in $M$.
\end{itemize}
\end{lemma}

We omit the easy proof of the next observation.

\begin{lemma}
\label{clonal-uniform}
The ground set of a matroid $M$ consists of a single
clonal class if and only if $M$ is uniform. Moreover, if
$M$ is uniform and $r(M),r^*(M)>0$,
then each element of $M$ has freedom $r(M)$
and cofreedom $r(M^*)$.
\end{lemma}

Elements $a$ and $b$ of a matroid $M$ are {\em incomparable}
if $a$ is not freer than $b$ and $b$ is not freer than $a$.
Note that an element 
$a$ has freedom 0 if and only if it is a loop of
$M$ and has cofreedom 0 if and only if it is a coloop of $M$.

\begin{lemma}
\label{bound-freedom}
Let $a$ be an element of the matroid $M$ that is neither a loop
nor a coloop. If the freedom
of $a$ is $\gamma$ and the cofreedom of $a$ is $\delta$,
then $M$ has a $U_{\gamma,\gamma+\delta}$-minor.
\end{lemma}

\begin{proof}
Note that $\gamma,\delta>0$.
Assume that $M$ is not uniform. By Lemma~\ref{clonal-uniform}, 
$M$ has an element
$b$ such that $a$ and $b$ are not clones.
Then either (i) $a$ and $b$ are incomparable, (ii) $a$ is strictly
freer than $b$, or (iii) $b$ is strictly freer than $a$.

Assume that either (i) or (ii) holds and consider the matroid
$M\ba b$. Then, by Lemma~\ref{minor-freedom}(i), $a$ has 
freedom at least $\gamma$ in $M\ba b$ and by 
Lemma~\ref{cominor-freedom}(ii), $a$ has cofreedom at least
$\delta$ in $M\ba b$. As $\gamma,\delta>0$, we see that
$a$ is neither a loop nor a coloop of $M\ba b$.
On the other hand, if (iii) holds, then we 
apply Lemma~\ref{minor-freedom}(ii) and
Lemma~\ref{cominor-freedom}(i) to obtain the same conclusion
for $M/b$.

Iterating the above procedure we eventually arrive at a 
minor $N$ of $M$ using $a$ whose ground set consists
of a single clonal class, where the freedom of 
$a$ in $N$ is $\gamma'\geq \gamma$ and the cofreedom is
$\delta'\geq \delta$. By Lemma~\ref{clonal-uniform}
$N$ is uniform, indeed $N\cong U_{\gamma',\delta'+\gamma'}$.
The lemma now follows from the fact that 
$U_{\gamma',\delta'+\gamma'}$ has a $U_{\gamma,\gamma+\delta}$-minor.
\end{proof}

\begin{corollary}
\label{bound-cor}
Let $q$ be a prime power and let $M$ be a matroid that
is either $GF(q)$-representable or is an excluded minor
for $GF(q)$-representability. 
\begin{itemize}
\item[(i)] If $e$ is not fixed in $M$, 
then $e$ has cofreedom at most $q$.
\item[(ii)] If $e$ is not cofixed in $M$, 
then $e$ has freedom at most $q$.
\end{itemize}
\end{corollary}

\begin{proof}
Assume that $e$ is not fixed in $M$. If 
$e$ is a coloop, then $e$ has cofreedom 0, so that (i)
holds. Assume that $e$ is not a coloop of $M$.
Then $e$ has freedom $k\geq 2$. Assume that $e$ has cofreedom
$q'>q$. By Lemma~\ref{bound-freedom},
$M$ has a $U_{k,q'+k}$-minor and hence a $U_{2,q+3}$-minor.
But this matroid is neither $GF(q)$-representable nor an 
excluded minor for $GF(q)$ contradicting the choice of
$M$. Thus (i) holds in this case too.
Part (ii) is the dual of (i).
\end{proof}

Recall that the {\em connectivity function} $\lambda_M$
of a matroid $M$ on $E$ is defined, for all subsets $X$ of $E$
by $\lambda_M(X)=r(X)+r(E-X)-r(M)$. We say that a partition
$(X,E-X)$ of $E$ is a $k$-{\em separation} if $\lambda_M(X)< k$.
If $\lambda_M(X)=k-1$, then the $k$-separation is {\em exact}. 
Recall also that the {\em coclosure} operator of $M$,
denoted $\cl^*$ is the closure operator of $M^*$. Thus
$x\in  \cl^*(A)$ if and only if $x\in\cl_{M^*}(A)$.

\begin{lemma}
\label{guts-coguts-freedom}
Let $e$ be an element of the matroid $M$ and $(A,B)$ be a 
$(t+1)$-separation.
\begin{itemize}
\item[(i)] If $e\in\cl(A-\{e\})\cap\cl(B-\{e\})$, then $e$ has 
freedom at most $t$ in $M$.
\item[(ii)] If $e\in\cl^*(A-\{e\})\cap \cl^*(B-\{e\})$, 
then $e$ has cofreedom at most $t$ in $M$.
\end{itemize}
\end{lemma}

\begin{proof}
Consider (i). Let $M'$ be an extension of $M$ by a set $Z$ such that
the members of $Z\cup\{e\}$ are clones. Then 
$Z\cup \{e\}$ is contained in $\cl_{M'}(A)\cap \cl_{M'}(B)$.
It now follows from the submodularity of the rank function that
$r_{M'}(Z\cup \{e\})\leq t$ so that $e$ has freedom at most
$t$ in $M$. Part (ii) is the dual of (i).
\end{proof}

The next lemma is \cite[Theorem~4.1]{govw}. As the result is vital
we repeat the short proof here.

\begin{lemma}
\label{vital}
Let $e$ be an element of a rank-$r$ matroid $M$. Let
$R$ be a $GF(q)$-representation of $M\ba e$
considered as a restriction of $PG(r-1,q)$.
Let $K$ be a flat of $PG(r-1,q)$ such that,
for each flat $F$ of $M$ in which $e$ is not a coloop,
$e\in F$ if and only if the flat of $PG(r-1,q)$ that is 
spanned by $F-\{e\}$ contains $K$. Then the rank of 
$K$ is at most the freedom of $e$ in $M$.
\end{lemma}

\begin{proof}
Let $\mathbb F$ be an infinite extension field of $GF(q)$
and let $\mathcal P$ be the projective space 
of rank $r$ over $\mathbb F$. Thus $\mathcal P$ contains
$PG(r-1,q)$. Let $K'$ be the flat of $\mathcal P$
that is spanned by $K$. Then, for each flat $F$ of $M$
in which $e$ is not a coloop, $e$ is in $F$ if and only if 
the flat of $\mathcal P$ that is spanned by $F-\{e\}$ contains
$K'$. Let $K^*$ denote the set of points $x$ of $K'$ for which
$R\cup \{x\}$ is an $\mathbb F$-representation of $M$.
Note that an element $x$ of $K'$ is in $K^*$ if and only if,
for each flat $F$ of $M$ not containing $e$, the point 
$x$ is not contained in the flat of $\mathcal P$ spanned
by $F$. There is a finite number of flats of $M$ that 
do not contain $e$. Therefore, by a simple comparison of
measures, $K^*$ spans $K'$. It is now straightforward to 
deduce that $K^*$ is spanned by some independent set $S$ such that
$S$  is a clonal class of the matroid $M'$ that is
represented by $R\cup S$. Note that $M'$ is an extension
of $M$ so that $|S|$ is at most the freedom of $e$ in $M$.
\end{proof}

\section{A certificate for a short proof of non-representability}

\begin{theorem}
\label{sufficient}
Let $q$ be a prime power and let
${\mathcal C}$ be a class of matroids 
such that the following hold.
\begin{itemize}
\item[(i)] For all nonempty matroids $M$ in $\mathcal C$,
there exists $e\in E(M)$ such that either $M\ba e$
or $M/e$ is in $\mathcal C$.
\item[(ii)] There exists an integer $t$ such that, if 
$M\in {\mathcal C}$ and $e\in E(M)$, then the following hold.
\begin{itemize}
\item[(a)] If $M/e\notin {\mathcal C}$,
then there is a $t$-separation $(A,B)$ in $M$ such that
$e\in\cl_M(A-\{e\})\cap \cl_M(B-\{e\})$.
\item[(b)] If $M\ba e\notin {\mathcal C}$, then there
is a $t$-separation $(A,B)$ in $M$ such that 
$e\in \cl^*_M(A-\{e\})\cap \cl^*_M(B-\{e\})$.
\end{itemize}
\item[(iii)] There exists an integer $s$ such that each matroid
$M\in {\mathcal C}$ has at most $s$ inequivalent representations
over $GF(q)$.
\item[(iv)] Each excluded minor for the class of 
$GF(q)$-representable matroids belongs to $\mathcal C$.
\end{itemize}
Then proving that an $n$-element matroid is not representable
over $GF(q)$ requires at most $O(n^2)$ rank evaluations.
\end{theorem}

\begin{proof}
Let $M$ be a 
matroid that is not representable over $GF(q)$.
In what follows, suppose that we have a {\em Claimant\ }
whose brief is to succinctly prove to an {\em Adjudicator\ } 
that $M$ is not $GF(q)$-representable.
The Claimant knows everything about $M$ but can only
reveal quadratically many rank-values to 
the Adjudicator.  The Claimant can find a minimal
minor $M'=M\ba D/C$ of $M$ that is not $GF(q)$-representable.
Now, for any $X\subseteq E(M')$, we have
$r_{M'}(X) = r_M(X\cup C) - r_M(C)$; thus
one rank evaluation for $N$ requires only two rank evaluations for $M$
(and if we need to make multiple rank evaluations for $M'$,
we only need to compute $r_M(C)$ once).
The Adjudicator concedes that it suffices to show that 
$N$ is not GF$(q)$-representable.
Thus we lose no generality in assuming that $M$ is an excluded
minor for $GF(q)$-representability.

Let $l=\max\{t,q\}$.

\begin{sublemma}
\label{sub}
If $N$ is a nonempty member of $\mathcal C$ and $N$ is either
an excluded minor for $GF(q)$-representability or is 
$GF(q)$-representable, then there is an element $e$ of
$N$ such that either:
\begin{itemize}
\item[(a)] $N\ba e\in\mathcal C$ and $e$ has freedom at most $l$, or
\item[(b)] $N/e\in\mathcal C$ and $e$ has cofreedom at most $l$.
\end{itemize}
\end{sublemma}

\begin{proof}
By property (i) of $\mathcal C$, there is an element $e\in E(N)$
such that either $N\ba e$ or $N/e$ is in $\mathcal C$. Up to duality
we may assume that $N\ba e$ is in $\mathcal C$. If the freedom 
of $e$ in $N$ is at most $l$, then the claim holds. Assume otherwise.
By Lemma~\ref{bound-cor}, $e$ is cofixed in $N$. Assume that
$N/e\notin\mathcal C$. Then, by property (ii) of $\mathcal C$,
there is a $t$-separation $(A,B)$ in $N$ such that 
$e\in\cl(A-\{e\})\cap \cl(B-\{e\})$. Then, by 
Lemma~\ref{guts-coguts-freedom}, $e$ has freedom at most $t\leq l$.
This contradiction shows that $N/e\in\mathcal C$. As 
$e$ is cofixed it has cofreedom $1<l$. Thus (b) is satisfied by $e$.
\end{proof}

By property (iv) of ${\mathcal C}$
our excluded minor $M$ belongs to ${\mathcal C}$. Say $|E(M)|=k$.
By \ref{sub} the Claimant can find a sequence
$M_0,M_1,\ldots, M_k=M$ of matroids in $\mathcal C$ such that
$M_0$ is empty, and
for each $i\in\{1,2,\ldots,k\}$ there is an element $e_i\in M_i$
such that either
\begin{itemize}
\item[] $M_{i-1}=M_i\ba e_i$ and $e_i$ has freedom
at most $l$ in $M_i$, or 
\item[] $M_{i-1}=M_i/e_i$ and $e_i$ has cofreedom at most
$l$ in $M_i$.
\end{itemize} 
For each $i\in \{1,2,\ldots,k\}$, let ${\mathcal R}_i$ be a 
complete set of inequivalent
$GF(q)$-representations of $M_i$; that is, any 
$GF(q)$-representation
of $M_i$ is equivalent to some representation in ${\mathcal R}_i$,
but no two representations in ${\mathcal R}_i$ are equivalent.
By property (iii) of $\mathcal C$, we have
$|{\mathcal R}_i|\leq s$ for all $i\in\{1,2,\ldots,k\}$.  
Moreover, since $M$ is not
$GF(q)$-representable, ${\mathcal R}_k$ is empty.  
The Claimant, who knows
everything about $M$, can determine 
$({\mathcal R}_1,{\mathcal R}_2,\ldots,{\mathcal R}_k)$.
The Claimants proof will consist of the sets
$({\mathcal R}_1,{\mathcal R}_2,\ldots,{\mathcal R}_k)$ 
along with a recursive argument
that each representation of $M_i$ is equivalent to
one in ${\mathcal R}_i$.

Suppose that the Adjudicator is already satisfied that
each $GF(q)$-representation of $M_{i-1}$ is equivalent to
some representation in ${\mathcal R}_{i-1}$. By duality
we may assume that $M_{i-1}=M_i\ba e_i$. 
Consider a representation $R_i$ of $M_i$. 
Evidently this representation is obtained by extending
a representation $R$ of $M_{i-1}$.
For each such representation we need to identify all
possible points $P$ in $PG(r-1,q)$ such that $R\cup\{p\}$
represents $M_i$ for $p\in P$. Suppose that
$S\subseteq E(M_i)-\{e_i\}$ with $e_i\in\cl_{M_i}(S)$.
Then the set $P$ is in the subspace of $PG(r-1,q)$
spanned by $S$. The Claimant will try to identify $P$ by 
considering the intersection of all such flats.
This is done inductively.
The Claimant constructs a sequence $K_0,K_1,\ldots,K_m$
of subspaces of PG$(r-1,q)$ as follows. Let $K_0 = PG(r-1,q)$.
For the flat $K_j$ one of the following holds.
\begin{itemize}
\item[1.] There is a set $S_j\subseteq E(M_i)-\{e\}$
and an element $a_j$ of $K_j$ 
such that $e_i$ is in the closure of $S_j$ in $M_i$
and $a_j$ is not spanned by $S_j$ in PG$(r-1,q)$.
In this case, the Claimant defines $K_{j+1}$ to be the intersection
of $K_j$ with the flat of PG$(r-1,q)$ spanned by $S_j$. 
\item[2.] For each  flat $F$ of $M_i$
containing $e_i$ such that $e_i$ 
is not a coloop of $M_i|F$, the flat $K_j$ is 
contained in the flat of PG$(r-1,q)$
that is spanned by $F-\{e\}$. Then $j=m$. 
\end{itemize}

Note that $m\le r$ and $K_m$ contains the set $P$ (which may
be empty).
The Claimant reveals the sets $(S_0,S_1,\ldots,S_{m-1})$ to the 
Adjudicator. Then, by revealing $O(r)$ rank values the Claimant
convinces the Adjudicator that
$e_i$ is in the closure of each of  $S_0,S_1,\ldots,S_{m-1}$. 
Given $S_0,S_1,\ldots,S_{m-1}$, the Adjudicator can then  determine
$K_0,K_1,\ldots,K_m$ efficiently using routine linear algebra.
Now we are in one of the following cases.
\begin{itemize}
\item[Case 1.] There is a set $S\subseteq E(M_i)-\{e_i\}$
such that $e_i$ is not in the closure of $S$ in $M_i$
but $K_m$ is contained in the flat spanned by $S$ in PG$(r-1,q)$.
\item[Case 2.] For each  flat $F$ of $M_i$ that does
not contain $e_i$, the flat $K_m$ is not contained in the flat
of $PG(r-1,q)$ that is spanned by $F$.
\end{itemize}

In Case 1 the Claimant can easily convince the Adjudicator
$R$ cannot be extended to a representation of $M_i$.
Indeed, two rank-values satisfy the 
Adjudicator that $e\in\cl_{M_i}(S)$ 
and the
Adjudicator can check that $K_m$ is spanned by $S$.
Now consider the second case. 
As $e_i$ has cofreedom at most $l$, it follows from Lemma~\ref{vital}
that $K_m$ has rank at most $l$. If $K_m$ has rank 0, then the 
Adjudicator concedes that $R$ cannot be extended to a representation
of $M_k$. Suppose that $K_m$ is nonempty. Then $K_m$
has at most $\frac{q^l-1}{q-1}$ elements. For an element
of $K_m$ that does not extend $R$ to a representation of 
$M_i$, the Claimant can reveal a single rank value to 
expose the fault, that is, the elements of 
$K_m$ that are not in $P$ can be identified with $O(1)$
rank evaluations.

As ${\mathcal R}_{i-1}$ has at most $s$ members, we need only
$O(r)$ rank evaluations to determine ${\mathcal R}_i$
from ${\mathcal R}_{i-1}$.
Therefore $O(|E(M)|^2)$ rank evaluations suffice
to prove that $M$ is not $GF(q)$-representable.

\end{proof}

\section{A short proof of non-representability for prime fields}

Let $k\geq 5$ be an integer. Then the definition of what it means
for a matroid to be $k$-{\em coherent} is given in \cite{gewh1}.
In fact $k$-coherence is a connectivity condition intermediate between
3-connectivity and 4-connectivity; $k$-coherent matroids are allowed
to have 3-separations but only in a controlled way.
The full definition takes some preparation and we will not give it
here. As a slight weakening of $k$-coherence we say that
a matroid is {\em near} $k$-coherent if it 
is connected and either $\si(M)$ or $\co(M)$ is 
$k$-coherent. The next theorem is \cite[Corollary~12.6]{gewh1}.

\begin{theorem}
\label{k-coherent}
Let $p$ be a prime number and $k\geq 5$ be an integer. Then the
following hold.
\begin{itemize}
\item[(i)] If $M$ is a nonempty near $k$-coherent matroid,
then there is an element $e\in E(M)$ such that either 
$M\ba e$ or $M/e$ is near $k$-coherent.
\item[(ii)] Let $M$ be a near $k$-coherent and $e\in E(M)$.
Then the following hold.
\begin{itemize}
\item[(a)] If $M/e$ is not near $k$-coherent, then there is a 
$4$-separation $(A,B)$ in $M$ such that 
$e\in\cl_M(A-\{e\})\cap\cl_M(B-\{e\})$.
\item[(b)] If $M\ba e$ is not near $k$-coherent, then there
is a $4$-separation $(A,B)$ in $M$ such that
$e\in\cl^*_M(A-\{e\})\cap \cl^*_M(B-\{e\})$.
\end{itemize}
\item[(iii)] There is an integer $\mu_p$ such that a 
near $k$-coherent matroid has at most $\mu_p$ inequivalent
representations over $GF(p)$.
\end{itemize}
\end{theorem}

The next theorem is a special case of a theorem of 
Ben David and Geelen \cite{dage}.
A {\em nested sequence} of $3$-separations of {\em length} $m$ 
in a matroid 
$M$ is a sequence $A_1,A_2,\ldots,A_m=E(M)$ of subsets of $E(M)$
such that, for all $i\in \{1,2,\ldots,m-1\}$, we have
$A_i$ is a proper subset of $A_{i+1}$ and $\lambda(A_i)\leq 2$. 

\begin{theorem}
\label{ben}
Let $p$ be a prime. Then there is an integer
$\nu_p$ such that an excluded minor for $GF(p)$-representability has
no nested sequence of $3$-separations of length $\nu_p$.
\end{theorem}

It is an immediate consequence of the definition of $k$-coherent
that a 3-connected matroid with no nested 
sequence of $3$-separations of
length $k$ is $k$-coherent. From
this fact, Theorem~\ref{ben}, and the fact that
$k$-coherent matroids are near $k$-coherent we
obtain the next corollary.

\begin{corollary}
\label{ben-cor}
Let $p$ be a prime number. Then an excluded minor for 
$GF(p)$-representability is $\nu_p$-coherent.
\end{corollary}

We now have all the ingredients for the proof of Theorem~\ref{certify}.
Let $p$ be a prime number and let $\mathcal C$ be the class
of near $\nu_p$-coherent matroids. Consider conditions (i)--(iv) of
Theorem~\ref{sufficient}. Condition (i) holds by 
Theorem~\ref{k-coherent}(i). Condition (ii) holds with $t=4$
by Theorem~\ref{k-coherent}(ii). Condition (iii) holds with
$s=\mu_p$ by Theorem~\ref{k-coherent}(iii). Finally condition
(iv) holds by Corollary~\ref{ben-cor}. Theorem~\ref{certify}
now follows from  Theorem~\ref{sufficient}.

\end{document}